\documentclass[a4paper,12pt]{book}
\usepackage{amsmath}
\usepackage[latin1]{inputenc}
\usepackage{amssymb}
\usepackage{graphicx}
\usepackage{makeidx}
\usepackage{mathrsfs}
\usepackage{amsfonts}
\usepackage[mathscr]{eucal}
\usepackage{amsmath}
\font\fmtitle=cmbx12 scaled \magstep2
\font\text=cmr10 at 12 truept
 at 11 truept
 at 10 truept

\begin{document}
 \thispagestyle{empty}
 
%\setcounter{page}{269}   \thispagestyle{empty}
%%%%%%% making running heads %%%%%%%%%%%%%%%%%%%%%%%%%%%
 \markboth
 {\it \centerline {Fractional Calculus in Linear Viscoelasticity}}
 {\it \centerline{Francesco MAINARDI}}  %% \hspace*{-20pt}}
%%%%%%% making running heads on even and odd pages %%%%%%%
%%%%%%%%%%% GENERAL DEFINITIONS
\def\pni{\par \noindent}
\def\vsh{\vskip 0.25truecm\noindent}
\def\vs{\vskip 0.5truecm}
\def\vvs{\vskip 1.0truecm}
\def\vvvs{\vskip 1.5truecm}
\def\vsp{\vsh\par}
\def\vsn{\vsh\pni}
\def\cen{\centerline}
\def\ra{\item{(a)\ }} \def\rb{\item{(b)\ }}   \def\rc{\item{(c)\ }}
\def\eg{{e.g.}\ } \def\ie{{i.e.}\ }
\def\sg{\hbox{sign}\,}
\def\sgn{\hbox{sign}\,}
\def\sign{\hbox{sign}\,}
\def\e{\hbox{e}}
\def\exp{\hbox{exp}}
%%%% MATHEMATICS
\def\ds{\displaystyle}
\def\dis{\displaystyle}
\def\q{\quad}    \def\qq{\qquad}
\def\lan{\langle}\def\ran{\rangle}
\def\l{\left} \def\r{\right}
\def\lt{\left} \def\rt{\right}
\def\lra{\Longleftrightarrow}
\def\arg{{\rm arg}}
\def\argz{{\rm arg}\, z}
\def\argG{{x^2/ (4\,a\, t)}}
\def\d{\partial}
 \def\dr{\partial r}  \def\dt{\partial t}
\def\dx{\partial x}   \def\dy{\partial y}  \def\dz{\partial z}
\def\rec#1{\frac{1}{#1}}
\def\log{{\rm log}\,}
\def\erf{{\rm erf}\,}     \def\erfc{{\rm erfc}\,}

%%%%%%%% SETS of NATURAL, REAL, COMPLEX NUMBERS : \NN, \RR, \CC
\def\NN{{\rm I\hskip-2pt N}}
\def\MM{{\rm I\hskip-2pt M}}
\def\RR{\vbox {\hbox to 8.9pt {I\hskip-2.1pt R\hfil}}\;}
\def\CC{{\rm C\hskip-4.8pt \vrule height 6pt width 12000sp\hskip 5pt}}
%%% IDENTITY
\def\II{{\rm I\hskip-2pt I}}
%%%%%%%%%%%%%%%%%%%%%%%%%%
%% DEFINITIONS of ERROR and EXPONENTIAL FUNCTIONS %%%%%
\def\erf{{\rm erf}\,}   \def\erfc{{\rm erfc}\,}
\def\exp{{\rm exp}\,} \def\e{{\rm e}}
\def\ss{{s}^{1/2}}   %% for LAPLACE TRANSFORMS
%%%%%%%%%%%%%%%%%%%%%%%%%%%%%%%%%%%%%%%%%%%%%%%%%%%%%%%
%%%%%%%%%%%%%%%%%%%%%%%%%%
\def\N{\bar N}  %%%%%%%%%%%%%
\def\ss{{s}^{1/2}} %%%%%%%
%%%%%%%%%%
%%%%%%%%%%%%%%%%%%%%%
%\def\sss{{s}^{1/2}}   %% for LAPLACE TRANSFORMS
\def\stt{{\sqrt t}}
\def\lst{{\lambda \,\stt}}
\def\Et{{E_{1/2}(\lst)}}
\def\u{\widetilde{u}}
\def\ul{\widetilde{u}} %%% Laplace Transform  (LT)
\def\uf{\widehat{u}} %%% Fourier Transform  (FT)
\def\bar{\widetilde}
\def\A{{\mathcal {A}}}
%%  justaposition symbols for Fourier, Laplace, Mellin transforms
\def\L{{\cal L}} %%% Laplace Transform !!!!
\def\F{{\cal F}} %%% Fourier Transform !!!!
\def\M{{\cal M}}  %%% Mellin Transform
\def\Fdiv{\,\stackrel{{\cal F}} {\leftrightarrow}\,}
  \def\Ldiv{\,\stackrel{{\cal L}} {\leftrightarrow}\,}
  \def\Mdiv{\,\stackrel{{\cal M}} {\leftrightarrow}\,}
%%%%%%%%%%% GREEN FUNCTIONS
\def\args{(x/ \sqrt{a})\, s^{1/2}}
\def\argsa{(x/ \sqrt{a})\, s^{\beta}}
\def\barr{\widetilde}
\def\argG{ x^2/ (4\,a\, t)}
\def\G{{\cal {G}}}
\def\Gc{{\cal {G}}_c}	\def\Gcs{\barr{\Gc}} %% CAUCHY PROBLEM
\def\Gs{{\cal {G}}_s}	\def\Gss{\barr{\Gs}} %% SIGNALLING PROBLEM
%%%%%% LAPLACE TRANSFORM
\def\f{\bar{f}}
\def\g{\bar{g}}
\def\u{\bar{u}}
%%%%%%%%%%%%%%%%%%
\centerline{{\fmtitle An Historical Perspective on}}
\vskip 0.10truecm  
\centerline{{\fmtitle Fractional Calculus in Linear Viscoelasticity}%%
\footnote{This paper has been published as a {\it Short Survey} in {\it Fract.Calc. Appl. Anal.}
{\bf 15} No 4 (2012),  712--717.}}
. 
%\vskip 0.10truecm
%% \text
\begin{center}
{\bf Francesco MAINARDI}
\vskip 0.10truecm \noindent
% \\
Deprtment of Physics, University of Bologna, and INFN, I-40126 Bologna, Italy
\\ E-mail: {\tt francesco.mainardi@bo.infn.it} 
\\ URL: {\tt http://www.fracalmo.org/mainardi/}
\end{center}
%%%%%%%%%%%%%%%%% SUMMARY %%%%%%%%%
%def\indice{\leaders\hbox to 1 em {\hss.\hss}\hfill}
%\def\hb#1{\hbox to 0.9 truecm{ \hss#1}}
\vskip 0.25truecm
\centerline{{\bf Abstract}}
\vskip 0.25truecm
\noindent 
The article provides an historical survey of the
early contributions on the applications of fractional calculus 
in linear viscoelasticty.
The period under examination covers four decades, since 1930's up to 1970's, 
and authors are from both Western and Eastern countries.
References to more recent contributions may be found in the  bibliography 
of the author's  recent book.
\\
This paper reproduces,  with Publisher's permission,  Section 3.5 of the book: 
F. Mainardi, {\it Fractional Calculus and Waves in Linear Viscoelasticity}, 
Imperial College Press -  London and World Scientific - Singapore, 2010 
\vskip .25truecm \noindent 
{\it AMS Subject Classification:}  
%% 01-XX, %% History of Mathematics
26A33,  %% Fractional Derivatives and Integrals  
%% 33E12 Mittag-Leffler functions
74-00, %% General reference works (handbooks, dictionaries, bibliographies, etc.)
%% 74-XX   %% mechanics of deformable bodies
74D05. %%  Linear constitutive equations
\vskip 0.25truecm
\section*{1. The first generation of pioneers of fractional calculus in viscoelasticity}
Linear viscoelasticity is certainly the field of the most extensive applications of fractional calculus, 
in view of its ability to model hereditary phenomena with long memory.
During the twentieth century a number of authors have (implicitly or explicitly) used the fractional calculus 
as an empirical method of describing the properties of viscoelastic materials. 

In the first half of that century the early contributors were: 
Gemant in USA, see [Gemant (1936); (1938)], 
Scott-Blair in England, see [Scott-Blair (1944); (1947); (1949)], 
Gerasimov and Rabotnov in the former Soviet Union, see [Gerasimov (1948)], [Rabotnov (1948)].

In 1950 Gemant published a series of 16 articles entitled Frictional Phenomena in 
Journal of Applied Physics since 1941 to 1943, which were collected in a book of the same title 
[Gemant (1950)]. 
In his eighth chapter-paper [Gemant (1942)], p. 220, he referred to his previous articles 
[Gemant (1936); (1938))] for justifying the necessity of fractional differential operators to 
compute the shape of relaxation curves for some elasto-viscous fluids. 
Thus, the words fractional and frictional were coupled, presumably for the first time, by Gemant. 
%%%%%%%%%%%%
% \newpage
%%%%%%%%%

Scott-Blair used the fractional calculus approach to model the observations made by 
[Nutting (1921); (1943); (1946)] that the stress relaxation phenomenon could be described by 
fractional powers of time. He noted that time derivatives of fractional order would simultaneously model 
the observations of Nutting on stress relaxation and those of Gemant on frequency dependence. 
It is quite instructive to cite some words by Scott-Blair quoted in [Stiassnie (1979)]: 
{\it I was working on the assessing of firmness of various materials 
(e.g. cheese and clay by experts handling them) these systems are of course both elastic 
and viscous but I felt sure that judgments were made not on an addition of elastic 
and viscous parts but on something in between the two so 
I introduced fractional differentials of strain with respect to time}. 
Later, in the same letter Scott-Blair added: 
{\it I gave up the work eventually, mainly because 
I could not find a definition of a fractional differential that would satisfy the mathematicians.}

The 1948 papers by Gerasimov and Rabotnov were published in Russian, 
so their contents remained unknown to the majority of western scientists up to the translation 
into English of the treatises by Rabotnov, see [Rabotnov (1969); (1980)]. 
Whereas Gerasimov explicitly used a fractional derivative to define his model of viscoelasticity 
(akin to the Scott-Blair model), 
Rabotnov preferred to use the Volterra integral operators with weakly singular kernels 
that could be interpreted in terms of fractional integrals and derivatives. 
After the appearance of the books by Rabotnov it has became common to speak about 
Rabotnov's theory of hereditary solid mechanics. 
The relation between Rabotnov's theory and the models of fractional viscoelasticity 
has been briefly recalled in the recent paper [Rossikhin and Shitikova (2007)]. 
According to these Russian authors, Rabotnov could express his models in terms of the operators of 
the fractional calculus, but he considered these operators only as some mathematical abstraction.

\section*{2. The second generation of pioneers of fractional calculus in viscoelasticity}
In the late sixties, formerly Caputo, see [Caputo (1966); (1967); (1969)], 
then Caputo and Mainardi, see [Caputo and Mainardi (1971a); (1971b)], 
explicitly suggested that derivatives of fractional order (of Caputo type) could be successfully 
used to model the dissipation in seismology and in metallurgy. 

In this respect the present author likes to recall a correspondence carried out between himself 
(as a young post-doc student) and the Russian Academician Rabotnov, related to two courses on 
Rheology held at CISM (International Centre for Mechanical Sciences, Udine, Italy) in 1973 and 1974, 
where Rabotnov was an invited speaker but without participating, see [Rabotnov (1973); (1974)]. 
Rabotnov recognized the relevance of the review paper [Caputo and Mainardi (1971b)], 
writing in his unpublished 1974 CISM Lecture Notes: 
{\it That's way it was of great interest for me to know the paper of Caputo and Mainardi from 
the University of Bologna published in 1971. 
These authors have obtained similar results independently without knowing the 
corresponding Russian publications.....}
 Then he added: 
 {\it The paper of Caputo and Mainardi contains a lot of experimental data of different authors 
 in support of their theory. On the other hand a great number of experimental curves obtained by 
 Postnikov and his coworkers as also by foreign authors can be found in numerous papers of Shermergor and Meshkov.} 
 Unfortunately, the eminent Russian scientist did not cite the 1971 paper by Caputo and Mainardi 
 (presumably for reasons independently from his willing) in the Russian and English editions of his later book 
 [Rabotnov (1980)].

Nowadays, several articles (originally in Russian) by Shermergor, Meshkov and their associated 
researchers have been re-printed in English in Journal of Applied Mechanics and Technical Physics 
(English translation of Zhurnal Prikladnoi Mekhaniki i Tekhnicheskoi Fiziki), 
see e.g. [Shermergor (1966)], [Meshkov et al. (1966)], [Meshkov (1967)], 
[Meshkov and Rossikhin (1968)],[Meshkov (1970)], [Zelenev et al. (1970)], 
[Gonsovskii and Rossikhin (1973)], 
available at the URL: {\tt http://www.springerlink.com/}.  
On this respect we cite the recent review papers [Rossikhin (2010)], [Rossikhin and Shitikova (2010)] 
where the works of the Russian scientists on fractional viscoelasticity are examined.

The beginning of the modern applications of fractional calculus in linear viscoelasticity is 
generally attributed to the 1979 PhD thesis by Bagley (under supervision of Prof. Torvik), 
see [Bagley (1979)], followed by a number of relevant papers, e.g. 
[Bagley and Torvik (1979); (1983a); (1983b)] and [Torvik and Bagley (1984)]. 
However, for the sake of completeness, one would recall also the 1970 PhD thesis of Rossikhin 
under the supervision of Prof. Meshkov, see [Rossikhin (1970)], 
and the 1971 PhD thesis of the author under the supervision of Prof. Caputo, summarized in 
[Caputo and Mainardi (1971b)].
 
To date, applications of fractional calculus in linear and nonlinear viscoelasticity 
have been considered by a great and increasing number of authors to whom we have tried to refer 
in the huge (but not exhaustive) bibliography  of our recent book [Mainardi (2010)].
Furthermore, in our book the interested reader can find an accessible introduction to fractional calculus 
and related special functions.   

\section*{Acknowledgments}

The author is grateful to Imperial College Press to allow him  to re-print in this short survey paper 
  Section 3.5 of his book:  F. Mainardi, {\it Fractional Calculus and Waves in Linear Viscoelasticity}, 
Imperial College Press -  London and World Scientific - Singapore, 2010.
\section*{References}

\vsn
 Bagley,  R.L. (1979). 
 {\it Applications of Generalized Derivatives to Viscoelasticity}, 
 Ph. D. Dissertation, Air Force Institute of  Technology.
\vsn
Bagley, R.L. and  Torvik,  P.J. (1979).  
A generalized derivative model for an elastomer damper, 
{\it Shock Vib. Bull.}  {\bf 49}, 135--143.
\vsn
 Bagley, R.L. and  Torvik,  P.J. (1983a). 
 A theoretical basis for the application  of fractional calculus, 
 {\it J. Rheology}  {\bf 27}, 201--210.
\vsn
Bagley, R.L. and  Torvik,  P.J. (1983b). 
Fractional calculus - A different  approach to the finite element analysis of viscoelastically damped structures, 
{\it AIAA Journal}  {\bf 21}, 741-748.
\vsn
Caputo, M. (1966). 
Linear models of dissipation whose $Q$ is almost frequency independent,    
{\it Annali di Geofisica}  {\bf 19}, 383-393.
\vsn
Caputo, M. (1967).  
Linear models of dissipation whose $Q$ is almost frequency  independent,  Part II, 
{\it Geophys. J. R. Astr. Soc.} {\bf 13}, 529-539.
 [Reprinted in  {\it Fract. Calc.  Appl. Anal.}  {\bf 11}, 4--14 (2008)] 
\vsn
 Caputo, M. (1969). 
 {\it Elasticit\`a e Dissipazione}, Zanichelli, Bologna. 
\vsn
 Caputo, M. and  Mainardi,  F. (1971a). 
  A new dissipation model based on memory mechanism,  
  {\it Pure Appl. Geophys. (PAGEOPH)}  {\bf 91}, 134-147. 
  [Reprinted in {\it Fract. Calc.  Appl. Analys.} {\bf 10}, 309--324 (2007)]
\vsn
 Caputo, M. and  Mainardi, F. (1971b).  
 Linear models of dissipation in  anelastic solids, 
 {\it Rivista del Nuovo Cimento (Ser. II)} {\bf 1}, 161--198.
\vsn
Gemant, A. (1936). 
A method of analyzing experimental results obtained from elastiviscous bodies, 
{\bf Physics}   {\bf 7},  311--317.  
\vsn
Gemant, A. (1938). 
On fractional differentials, 
{\it Phil. Mag. (Ser. 7)}  {\bf 25}, 540--549.
 \vsn
Gemant, A. (1942). 
Frictional phenomena: VIII,  
{\it J. Appl. Phys.}  {\bf 13}, 210-221. 
\vsn
Gemant, A. (1950).  
{\it Frictional Phenomena}, Chemical Publ. Co, Brooklyn N.Y.
\vsn
 Gerasimov, A. (1948). 
 A generalization of linear laws of deformation and its applications to problems of internal friction, 
 {\it Prikl. Matem. i Mekh. (PMM)}  {\bf 12} No 3, 251--260. [in Russian]
\vsn
Gonsovski, V.L. and Rossikhin, Yu.A. (1973). 
Stress waves in a viscoelastic medium with a singular hereditary kernel, 
{\it J. Appl. Mech. Tech. Physics} {\bf 14} No 4, 595--597. 
\vsn
 Mainardi, F. (2010). 
 {\it Fractional Calculus and Waves in Linear Viscoelasticity},
  Imperial College Press, London. 347 pages
[ISBN 978-1-84816-329-4].
 \\  See  {\tt http://www.worldscientific.com/worldscibooks/10.1142/p614}
\vsn 
Meshkov, S.I., Pachevskaya, G.N.  and Shermergor, T.D. (1966).
Internal friction described with the aid of fractionally-exponential kernels, 
{\it J. Appl. Mech. Tech. Physics}  {\bf 7} No 3, 63--65.
\vsn
 Meshkov, S.I. (1967). 
 Description of internal friction in the memory theory of elasticity using kernels with a weak singularity,  
 {\it J. Appl. Mech. Tech. Physics}  {\bf 8} No 4, 100--102.
\vsn
Meshkov, S.I. and  Rossikhin, Yu.A. (1968). 
Propagation of acoustic waves in a hereditary elastic medium, 
{\it J. Appl. Mech. Tech. Physics} {\bf  9} No 5, 589--592.
\vsn
Meshkov, S.I. (1970). 
The integral representation of fractionally exponential functions 
and their application  to dynamic problems of linear viscoelasticity,  
{\it J. Appl. Mech. Tech. Physics}  {\bf 11} No 1, 100--107.
\vsn
Nutting,  P.G. (1921). 
A new general law of deformation, 
{\it J. Frankline Inst.} {\bf 191}, 679--685.
\vsn
 Nutting, P.G. (1943). 
 A  general stress-strain-time formula, 
 {\it J. Frankline Inst.}  {\bf 235}, 513--524.
\vsn
 Nutting, P.G. (1946). 
 Deformation in relation to time, pressure and temperature, 
 {\it J. Frankline Inst.}  {\it 242}, 449-458.
\vsn
 Rabotnov,  Yu.N. (1948). 
 Equilibrium of an elastic medium with after effect, 
 {\it Prikl. Matem. i Mekh. (PMM)}  {\bf 12} No 1, 81-91. [in Russian]
\vsn
 Rabotnov, Yu.N.  (1969).
 {\it Creep Problems in Structural Members},
  North-Holland, Amsterdam.
 [English translation of  the 1966 Russian edition]
\vsn
 Rabotnov, Yu.N. (1973). 
 {\it On the Use of Singular Operators in the Theory of Viscoelasticity}, Moscow, 1973, pp. 50. 
 Unpublished Lecture Notes  for the CISM course on Rheology   held in Udine, October 1973.
 [{\it http://www.cism.it}]
 \vsn
Rabotnov, Yu.N. (1974).  
{\it Experimental Evidence of the Principle of Hereditary in Mechanics  of Solids}, Moscow, 1974, pp. 80.  
Unpublished Lecture Notes  for the CISM course on  Experimental Methods in Mechanics, A) Rheology,   
held in Udine, 24-29 October 1974.
 [{\it http://www.cism.it}]
\vsn
Rabotnov,  Yu.N. (1980). 
{\it Elements of Hereditary Solid Mechanics},
  MIR, Moscow. [English translation, revised from the 1977 Russian]
\vsn
 Rossikhin, Yu.A. (1970).  
 {\it Dynamic Problems of Linear Viscoelasticity connected with the Investigation  of Retardation and 
 Relaxation Spectra},
PhD Dissertation, Voronezh Polytechnic Institute, Voronezh.[in Russian]
\vsn
Rossikhin, Yu.A. (2010).
Reflections on two parallel ways in the progress of fractional calculus in mechanics of solids, 
{\it Appl. Mech. Review}  {\bf 63},  010701/1--12.
\vsn
Rossikhin, Yu.A.  and  Shitikova, M.V.  (2007).
Comparative analysis of viscoelastic models involving fractional derivatives of different orders,  
{\it  Fract. Calc.  Appl. Anal.} {\bf 10} No 2,  111--121. %%  )No 2
\vsn
Rossikhin, Yu.A.  and  Shitikova, M.V. (2010). 
Applications of fractional calculus to dynamic problems of solid mechanics: novel trends and recent results, 
{\it Appl. Mech. Review}  {\bf 63},  010801/1--52.
\vsn  
Scott-Blair, G.W.  (1944). 
Analytical and integrative aspects  of the stress-strain-time problem,  
{\it J. Scientific Instruments}  {\bf 21}, 80--84.
\vsn
 Scott-Blair, G.W.  (1947). 
 The role of psychophysics in rheology, 
{\it J. Colloid Sci.}  {\bf 2}, 21-32. 
\vsn
Scott-Blair,  G.W. (1949).  
{\it Survey of General and Applied Rheology},
 Pitman, London.
\vsn %% \bibitem[Shermergor (1966)]{Shermergor PMTF66}
 Shermergor, T.D. (1966).
 On the use of fractional differentiation operators for the description of 
 elastic--after effect properties of materials,
 {\it J. Appl. Mech.  Tech. Phys.} {\bf 7} No 6, 85--87.
 % {\it Prikl. Mat. Tekh. Fiz.} {\bf 6}, 118-121.
%% DOI 10.1007/BF00914347 published on line 3 Jan 2005 in Springer Link
\vsn
 Stiassnie, M.  (1979). 
 On the application of fractional calculus on the formulation  of viscoelastic models, 
 {\it Appl. Math. Modelling}  {\bf 3}, 300--302.
\vsn %% \bibitem[Torvik and  Bagley (1984)]{Torvik-Bagley 84}
 Torvik, P.J.  and  Bagley, R.L. (1984).
On the appearance of the fractional derivatives in the behavior of
  real materials,
{\it ASME J. Appl. Mech.} {\bf 51}, 294--298.
\vsn %% \bibitem[Zelenev {\it et al}. (1970)]{Zelenev-Meshkov-Rossikhin 70}
Zelenev, V.M., Meshkov, S.I. and  Rossikhin, Yu.A. (1970).
Damped vibrations of hereditary - elastic systems with weakly singular kernels,
 {\it J. Appl. Mech. Tech. Physics} {\bf 11} No 2, 290--293.

\end{document}